\documentstyle{amsart} 
\input epsf.sty

\newcommand{\be}{\begin{enumerate}}

\newcommand{\ca}{{\cal A}}

\newcommand{\cd}{{\cal D}}

\newcommand{\ck}{{\cal K}}
\newcommand{\cl}{{\cal L}}

\newcommand{\cs}{{\cal S}}

\newcommand{\da}{\Delta}

\newcommand{\dc}{{\Bbb C}}
\newcommand{\dr}{{\Bbb R}}
\newcommand{\dq}{{\Bbb Q}}
\newcommand{\dz}{{\Bbb Z}}

\newcommand{\ee}{\end{enumerate}}

\newcommand{\U}{\text{\bf U}}
\newcommand{\GL}{\text{\bf GL}}

\newcommand{\lam}{\lambda}
\newcommand{\Lam}{\Lambda}

\newcommand{\lmu}{(\lam,\mu)}
\newcommand{\Lmu}{(\Lam,\mu)}

\newcommand{\pl}{\Pi_\lambda}
\newcommand{\pr}{\noindent{\bf Proof. }}

\newcommand{\ra}{\rightarrow}

\newcommand{\sig}{\sigma}
\newcommand{\sla}{\Sigma_\lam}
\newcommand{\sm}{\setminus}

\newcommand{\st}{\text{\rm St}}
\newcommand{\susp}{\text{\rm susp}}

\newcommand{\tb}{\tilde\beta}

\newcommand{\ti}{\tilde }

\newcommand{\wti}{\widetilde }
\newcommand{\wsl}{{\wti\Sigma}_\lam}

\newcommand{\xlm}{X_{\lambda,\mu}}
\newcommand{\xLm}{X_{\Lambda,\mu}}

\newtheorem{thm}{Theorem}[section]
\newtheorem{df}  [thm]{Definition}

\newtheorem{crl} [thm]{Corollary}
\newtheorem{prop}[thm]{Proposition}
\newtheorem{conj}[thm]{Conjecture}

\numberwithin{equation}{section}

\title [Rational homology of spaces of complex polynomials]
       {Rational homology of spaces of complex \\
       monic polynomials with multiple roots}
 
   \author{Dmitry N. Kozlov}
                
  \date{\today \\[0.05cm]
 \hskip15pt Research at IAS was supported by von Hoffmann, Arcana Foundation.}

\address{Institute for Advanced Study, Olden Lane, Princeton, NJ 08540, U.S.A.}

\email{kozlov@@math.ias.edu, kozlov@@math.kth.se} 
 
\begin{document}

\begin{abstract} 
   We study rational homology groups of one-point compactifications of spaces
of complex monic polynomials with multiple roots. These spaces are indexed by
number partitions. A~standard reformulation in terms of quotients
of orbit arrangements reduces the problem to studying certain triangulated 
spaces~$\xlm$. 

   We present a~combinatorial description of the cell structure of~$\xlm$ using 
the language of marked forests. As applications we obtain a~new proof of 
a~theorem of Arnold and a~counterexample to a~conjecture of Sundaram and Welker, 
along with a~few other smaller results.
\end{abstract}

\maketitle
              
              \section{Introduction.}
              
   Let $n$ be an~integer, $n\geq 2$. We view $n$-dimensional complex space $\dc^n$ 
as the space of all monic polynomials with complex coefficients of degree $n$ by 
identifying $a=(a_0,\dots,a_{n-1})\in\dc^n$ with $f_a(z)=z^n+a_{n-1}z^{n-1}+\dots+a_0$. 
To each $\lam=(\lam_1,\dots,\lam_t)\vdash n$ one can associate a~topological space 
as follows (we refer the reader to the subsection~\ref{s2.1} for a~description of our 
conventions on the terminology of number and set partitions).
\begin{df}
  $\wsl$ is the set of all $a\in\dc^n$, for which the roots of $f_a(z)$ can be 
partitioned into sets of sizes $\lam_1,\dots,\lam_t$, so that within each set 
the roots are equal. Clearly, $\wsl$ is a~closed subset of $\,\dc^n$. Let $\sla$ 
be the one-point compactification of~$\wsl$.
\end{df}

In this paper we shall focus on the (reduced) rational Betti numbers of spaces~$\sla$.
In~\cite{Ar2}, V.I.~Arnold has computed $\tb_*(\sla,\dq)$ for $\lam=(k^m,1^{n-km})$.
\begin{thm}\label{arn} {\rm (\cite{Ar2}).}
  Let $\lam=(k^m,1^{n-km})$ for some natural numbers $k\geq 2$, $m$, and $n\geq km$.
Then
\begin{equation}\label{are}
 \tb_i(\sla,\dq)=\begin{cases} 1,& \text{for } i=2l(\lam);\\
                               0,& \text{otherwise.}
  \end{cases}
\end{equation}
\end{thm}

  In \cite{SuW} Sundaram and Welker conjectured that 
   
\begin{conj}\label{c3}
  For any number partition $\lam$, $\tb_i(\sla,\dq)=0$ unless $i=2l(\lam)$.
\end{conj}

  In this paper we shall give a new, combinatorial proof of the Theorem~\ref{arn}
 and disprove Conjecture~\ref{c3}.

  To do that, we shall introduce a family of topological spaces $\xlm$,
indexed by pairs of number partitions $(\lam,\mu)$, satisfying $\lam\vdash\mu$. 
$\xlm$ will be defined so that the following equality is satisfied

\begin{equation}\label{rswf}
  \tb_i(\sla,\dq)=\sum_{\lam\vdash\mu\vdash n}\tb_{i-2l(\mu)-1}(\xlm,\dq).
\end{equation}
  
  Here is the summary of the paper.

 In Section 2 we introduce terminology of number and set partitions,
subspace arrangements and their intersection lattices, and order
complexes of posets.

 In Section 3 we define the topological spaces $\xlm$ and derive \eqref{rswf}.

 In~Section~\ref{s3} we give a~combinatorial description of the cell 
structure of the triangulated spaces~$\xlm$ in terms of marked forests, see
Theorem~\ref{main1}. This description is the backbone of the paper, it serves 
as both language and intuition for the material in the subsequent sections.
One consequence of Theorem~\ref{main1} is that homology groups of $\xlm$ may be
computed from a~chain complex, whose components are freely generated by marked
forests, and the boundary operator is described in terms of a~combinatorial
operation on such forests (deletion of level sets).

 In Section \ref{s5} we prove a~general theorem about collapsibility of 
certain triangulated spaces. The argument is along the lines of \cite{Ko2} which 
has been extended to cover the result of Arnold and clarified by using
the combinatorial cell description, 
from Section~\ref{s3}, of~$\xlm$. We would like to mention that for the special 
cases $\lam=(k,1^{n-k})$ and $\lam=(k^m)$ the Theorem~\ref{arn} was also reproved 
in~\cite{SuW} by using the Theorem~\ref{swt} (as Theorem~\ref{tgen} shows, the 
case $\lam=(k^m)$ is especially simple). However, the Theorem~\ref{main2} is 
the first combinatorial (modulo Theorem~\ref{swt}) proof of the result of Arnold 
in the general case.

 In Section~\ref{s4} we disprove the conjecture of Sundaram and Wel\-ker.
Besides giving a~counterexample we prove this conjecture for a~class of
number partitions, which we call generic partitions.

              \section{Terminology.}

\subsection{Number and set partitions.}\label{s2.1}$\,$
 
\noindent 
 Let $n$ be a~natural number. We denote the set $\{1,\dots,n\}$ by $[n]$. 
A~{\it number partition} of~$n$ is a~set $\{\lam_1,\dots,\lam_t\}$ of natural 
numbers, such that $\lam_1+\dots+\lam_t=n$. The usual con\-ven\-tion is to 
write $\lam=(\lam_1,\dots,\lam_t)$, where $\lam_1\geq\dots\geq\lam_t$, and 
$\lam\vdash n$. The {\it length} of $\lam$, denoted $l(\lam)$, is the number 
of components of $\lam$, say, in the previous sentence $l(\lam)=t$. We also 
use the power notation: $(n^{\alpha_n},\dots,1^{\alpha_1})=(\underbrace
{n,\dots,n}_{\alpha_n},\dots,\underbrace{1,\dots,1}_{\alpha_1})$. 

  For a~set~$S$, a~{\it set partition} $\pi$ is a~set $\{S_1,\dots,S_t\}$, 
where $S_i$'s are subsets of $S$, such that $\cup_{i=1}^t S_i=S$ and 
$S_i\cap S_j=\emptyset$ for $1\leq i<j\leq t$. We write $\pi=(S_1,\dots,S_t)$, 
where $|S_1|\geq\dots\geq|S_t|$, and $\pi\vdash S$. Observe that we do
not distinguish the set partitions differing only in the order of the sets, 
e.g., $(\{1,2\},\{3,4\})$ and $(\{3,4\},\{1,2\})$ are the same. Whenever we 
write $\pi\vdash[n]$, it implicitly implies that $\pi$ is a~set partition, 
as opposed to a~number partition. A~set partition $\pi\vdash[n]$, 
$\pi=(S_1,\dots,S_t)$, is said to have {\it type} $\lam$, where 
$\lam\vdash n$ is the~number partition $\lam=(|S_1|,\dots,|S_t|)$. 

For two set partitions $\pi,\ti\pi\vdash S$, $\pi=(S_1,\dots,S_t)$, $\ti\pi=
(\ti S_1,\dots,\ti S_q)$ we write $\pi\vdash\ti\pi$ if there exists 
$\iota\vdash[t]$, $\iota=\{I_1,\dots,I_q\}$, such that $\ti S_i=
\cup_{j\in I_i} S_j$ for $i\in[q]$. Analogously, for two number partitions 
$\lam=(\lam_1,\dots,\lam_t)$, $\mu=(\mu_1,\dots,\mu_q)$ we write $\lam\vdash\mu$ 
if there exists $\iota\vdash[t]$, $\iota=\{I_1,\dots,I_q\}$, such that $\mu_i=
\sum_{j\in I_i}\lam_j$ for $i\in[q]$. Clearly $\pi\vdash[n]$ and $\lam\vdash n$ 
are special cases of these notations. Finally observe that if $\pi,\ti\pi$ are 
two set partitions, such that $\pi\vdash\ti\pi$, then $(\text{type }\pi)\vdash
(\text{type }\ti\pi)$.

\subsection{Subspace arrangements and their intersection lattices.}
\begin{df}
  A set $\ca=\{\ck_1,\dots,\ck_t\}$ of linear complex subspaces of $\dc^n$, such that
$\ck_i\not\subseteq\ck_j$ for $i\neq j$, is called a~(central complex) {\bf subspace
arrangement} in~$\dc^n$. 
\end{df}

  The intersection data of a~subspace arrangement may be represented by a~poset.

\begin{df}
  To a~subspace arrangement $\ca$ in $\dc^n$ we associate a~partially ordered set
$\cl_\ca$, called the {\bf intersection lattice} of $\ca$. The set if elements of 
$\cl_\ca$ is $\{K\subseteq\dc^n\,|\,\exists\, I\subseteq[t],\text{ such that }
\bigcap_{i\in I} \ck_i=K\}\cup\{\dc^n\}$ with the order given by reversing inclusions: 
$x\leq_{\cl_\ca}y$ iff $x\supseteq y$. That is, the minimal element of $\cl_\ca$ 
is~$\dc^n$, also customarily denoted $\hat 0$, and the maximal element is 
$\bigcap_{K\in\ca}K$.
\end{df}  

  Let $V_\ca=\bigcup_{i=1}^t\ck_i$. If $V_\ca$ is invariant under the action 
of some finite group $G\subset\GL_n(\dc)$, then we say that $G$ acts on $\ca$. 
In that case, $\Gamma_\ca^G$ denotes the one-point compactification of~$V_\ca/G$. 
For $x\in\cl_\ca$, $\st_x\subseteq G$ denotes the stabilizer of~$x$.

\subsection{Order complexes of posets.}

\begin{df}
  For a~poset $P$, let $\da(P)$ denote the nerve of~$P$ viewed as a~category 
in the usual way: it is a~simplicial complex with $i$-dimensional simplices 
corresponding to chains of $i+1$ elements of $P$ (chains are totally ordered 
sets of elements of~$P$). In particular, vertices of~$\da(P)$ correspond to 
the elements of~$P$. We call $\da(P)$ the {\bf order complex} of~$P$.
\end{df}

  For an~arbitrary poset $P$ and $x,y\in P$, $x<y$, let $P(x,y)$ denote the 
subposet of~$P$ consisting of all $z\in P$, such that $x<z<y$.

              \section{Orbit arrangements and spaces $\xlm$.}
  
\subsection{Reformulation in the language of orbit arrangements.}$\,$

\noindent
  Following \cite{SuW} we shall give a~different interpretation of the numbers
$\tb_i(\sla,\dq)$, for general $\lam$. First, let us observe that the symmetric 
group~$\cs_n$ acts on~$\dc^n$ by permuting the coordinates, so we can consider 
the space $\dc^n/\cs_n$ endowed with the quotient topology. It is a~classical 
fact that the map $\phi:\dc^n\ra\dc^n/\cs_n$, mapping a~polynomial to the 
(unordered) set of its roots, is a~homeomorphism, which extends to the one-point 
compactifications. Therefore $\dc^n\cong\dc^n/\cs_n$ and 
$\sla\cong\phi(\sla)=\phi(\wsl)\cup\{\infty\}$.

   $\phi(\wsl)$ can be viewed as the configuration space of $n$ unmarked points 
on~$\dc$ such that the number partition given by the coincidences among the points 
is refined by~$\lam$. For example, $\phi(\wti\Sigma_{(2,1^{n-2})})$ is the 
configuration space of $n$ unmarked points on $\dc$ such that at least 2 points 
coincide. Using this point of view, $\phi(\sla)$ can be described in the language 
of orbit subspace arrangements.

\begin{df}
   For $\pi\vdash[n]$, $\pi=(S_1,\dots,S_t)$, $S_j=\{i_1^j,\dots,i_{|S_j|}^j\}$, 
$1\leq j\leq t$, $K_\pi$ is the~sub\-space given by the equations 
$x_{i_1^1}=\dots=x_{i_{|S_1|}^1},\dots,x_{i_1^t}=\dots=x_{i_{|S_t|}^t}$.
  For $\lam\vdash n$, set $I_\lam=\{\pi\vdash[n]\,|$ $\text{\rm type}\,(\pi)=\lam\}$
and define $\ca_\lam=\{K_\pi\,|\,\pi\in I_\lam\}$. $\ca_\lam$'s are called
{\bf orbit arrangements}.
\end{df}
   The orbit arrangements were introduced in \cite{Bj} and studied in further detail 
in~\cite{Ko1}. They provide the appropriate language to describe $\phi(\sla)$, indeed
\begin{equation}\label{fisla}
  \phi(\sla)=\Gamma_{\ca_\lam}^{\cs_n}.
\end{equation} 

   An~important special case is that of the braid arrangement 
$\ca_{n-1}=\ca_{(2,1^{n-2})}$, which corresponds under $\phi$ to 
${\wti\Sigma}_{(2,1^{n-2})}$, the~space of all monic complex polynomials of 
degree~$n$ with at least one multiple root. The name "braid arrangement" stems 
from the fact that $\dc^n\sm V_{\ca_{n-1}}$ is a~classifying space of the 
colored braid group, see~\cite{Ar1}. The intersection lattice $\cl_{\ca_{n-1}}$ 
is usually denoted $\Pi_n$. It is the~poset consisting of all set partitions 
of~$[n]$, where the partial order relation is refinement. Furthermore, for 
$\lam\vdash n$, the intersection lattice of $\ca_\lam$ is denoted $\Pi_\lam$. 
It is the~subposet of $\Pi_n$ consisting of all elements which are joins of 
elements of type~$\lam$, with the minimal element $\hat 0$ attached. 
          
\subsection{Applying Sundaram-Welker formula.}$\,$
  
\noindent   
 The following 
formula of S.~Sun\-da\-ram and V.~Welker, \cite{SuW}, is vital for our approach. 

\begin{thm}\label{swt}{\rm (\cite[Theorem 2.4(ii) and Lemma 2.7(ii)]{SuW}).}\newline 
Let $\ca$ be an~arbitrary subspace arrangement in $\dc^n$ with an~action of a~finite
group $G\subset\U_n(\dc)$. Let $\cd_\ca$ be the intersection of $V_\ca$ with the 
$(2n-1)$-sphere (often called the link of $\ca$). Then there is the following 
isomorphism of $G$-modules.
\begin{equation}\label{swf}
  \wti H_i(\cd_\ca,\dq)\cong_G\bigoplus_{x\in\cl_\ca^{>\hat 0}/G} 
  \text{\rm Ind}_{\st_x}^G (\wti H_{i-\dim x}(\da(\cl_\ca(\hat 0,x)),\dq)),
\end{equation}
  where the sum is taken over representatives of the orbits of $G$ in
$\cl_\ca\sm\{\hat 0\}$, under the action of $G$, one representative for each orbit.
\end{thm} 
   
   Clearly $\Gamma_\ca^G\cong\susp\,(\cd_\ca/G)$. Recall that if a~finite group~$G$ 
acts on a~finite cell complex $K$ then $\tb_i(K/G,\dq)$ is equal to the multiplicity 
of the trivial representation in the induced representation of $G$ on the $\dq$-vector 
space $\wti H_i(K,\dq)$, see for example \cite[Theorem 1]{Co}, \cite{Br}. Hence, it 
follows from \eqref{swf}, and the Frobenius reciprocity law, that 
\begin{equation}\label{star}
  \tb_i(\Gamma_\ca^G,\dq)=\sum_{x\in\cl_\ca^{>\hat 0}/G} 
  \tb_{i-\dim x-1}(\da(\cl_\ca(\hat 0,x))/{\st_x},\dq).
\end{equation}
   
\subsection{Spaces $\xlm$ and their properties.}$\,$

\noindent
  Let us now restate this identity in the special case of orbit arrangements. 
As mentioned above, the intersection lattice of $\ca_\lam$ is $\Pi_\lam$.
It has an~action of the symmetric group $\cs_n$, which, for any $\pi\in\Pi_\lam$ 
induces an~action of $\st_\pi$ on $\da(\Pi_\lam(\hat 0,\pi))$. 

{\bf Notation. }{\it Let $\xlm$ denote the topological space 
$\da(\pl(\hat 0,\pi))/\st_\pi$, where the set partition~$\pi$ has type $\mu$.
If there is no set partition $\pi\in\pl$ of type $\mu$, i.e., if $\mu$ 
cannot be obtained as a join of $\lambda$'s, then let $\xlm$ be a~point.}

For fixed $\mu$, the space $\xlm$ does not depend on the choice of 
$\pi$. Observe that $\xlm$ is in general not a~simplicial complex, however it is 
a~triangulated space, (a~regular CW complex with each cell being a~simplex, 
see~\cite[Chapter~I, Section~1]{GeM}), with its cell structure inherited from the 
simplicial complex $\da(\Pi_\lam(\hat 0,\pi))$. In general, whenever $G$ is 
a~finite group which acts on a~poset $P$ in an~order-preserving way, $\da(P)/G$ 
is a~triangulated space whose cells are orbits of simplices of $\da(P)$ under 
the action of $G$; this is obviously not true in general for an~action of 
a~finite group on a~finite simplicial complex.
        
  Clearly, \eqref{fisla} together with \eqref{star}, and the fact that
$\phi$ is a~homeomorphism, implies \eqref{rswf}.
Let us quickly analyze \eqref{rswf}. $X_{\lam,\lam}=\emptyset$ makes
a~contribution~1 in dimension $2l(\lam)$. Assume $\mu\neq\lam$, then $1\leq l(\mu)
\leq l(\lam)-1$ and $\xlm\neq\emptyset$. $\dim\xlm=l(\lam)-l(\mu)-1$, hence 
$\tb_{i-2l(\mu)-1}(\xlm,\dq)=0$ unless $0\leq i-2l(\mu)-1\leq l(\lam)-l(\mu)-1$,
that is $2l(\mu)+1\leq i\leq l(\lam)+l(\mu)$. It follows from \eqref{rswf} that
$$\tb_{2l(\lam)}(\sla,\dq)=1,\text{  and }\tb_i(\sla,\dq)=0\text{ unless }
3\leq i\leq 2l(\lam).$$

 The purpose of this paper is to investigate the values $\tb_i(\sla,\dq)$ for
$3\leq i\leq 2l(\lam)-1$, by studying $\tb_i(\xlm,\dq)$. We shall prove that 
the latter are equal to~0 for a~certain set of pairs $(\lam,\mu)$, 
$\lam\vdash\mu$, of partitions, including the case in Theorem~\ref{arn}, 
($\lam=(k^m,1^{n-km})$, $\mu$ is arbitrary such that $\lam\vdash\mu$),  
and we shall give an~example that this is not the case in general.

              \section{The cell structure of $\xlm$ and marked forests}
          \label{s3}
                 
 \subsection{The terminology of marked forests.}$\,$

\noindent  
  In order to index the simplices of $\xlm$ we need to introduce some terminology 
for certain types of trees with additional data. For an~arbitrary forest of rooted
trees $T$ (we only consider finite graphs), let $V(T)$ denote the set of the 
vertices of~$T$, $R(T)\subseteq V(T)$ denote the set of the roots of~$T$ and 
$L(T)\subseteq V(T)$ denote the set of the leaves of~$T$. For any integer $i\geq 0$,
let $l_i(T)$ be the number of $v\in V(T)$ such that, $v$ has distance $i$ to 
the root in its connected component. 

\begin{df}
  A forest of rooted trees $T$ is called a~{\bf graded forest of rank $r$} if 
$l_{r+2}(T)=0$, $l_{r+1}(T)=|L(T)|$, and the sequence $l_0(T),\dots,l_{r+1}(T)$ 
is strictly increasing.  
\end{df}

For $v,w\in V(T)$, $w$ is called {\it a~child} of $v$ if there is an~edge between $w$ 
and $v$ and the unique path from $w$ to the corresponding root passes through $v$.
For $v\in V(T)$, we call the distance from $v$ to the closest leaf the {\it height}
of $v$. For example, in a~graded forest of rank $r$, leaves have height 0 and roots
have height $r+1$.

\begin{df}\label{df2.2}
  A {\bf marked forest of rank $r$} is a~pair $(T,\eta)$, where $T$ is a~graded
forest of rank $r$ and $\eta$ is a~function from $V(T)$ to the set of natural 
numbers such that for any vertex $v\in V(T)\sm L(T)$ we have 
\begin{equation}\label{mc1}
\eta(v)=\sum_{w\in\text{\rm{children}}(v)}\eta(w).
\end{equation}
\end{df}

  We remark that the set of the marked forests of rank $r$, such that not all 
leaves have label~1, is equal to the set of graded forests of rank $r+1$. 
Indeed, instead of labeling the vertices with natural numbers so that~\eqref{mc1} 
is satisfied, one can as well attach a~new level of leaves so that each "old leaf" 
$v$ has $\eta(v)$ children. Then the old labels will correspond to the numbers of 
the new leaves below each vertex. For our context it is more convenient to use 
labels rather than auxiliary leaves, i.e.,~it is more handy to label all vertices 
rather than just the leaves, so we stick to the terminology of Definition~\ref{df2.2}.

   For a~marked forest $(T,\eta)$ of rank $r$ and $0\leq i\leq r+1$, we have 
a~number partition $\lam_i(T,\eta)=\{\eta(v)\,|\,v\text{ has height }i\}$. Clearly 
$\lam_0(T,\eta)\vdash\dots\vdash\lam_r(T,\eta)\vdash\lam_{r+1}(T,\eta)$.

\begin{df}\label{lmd}
   Let $\lam\vdash\mu\vdash n$, $\lam\neq\mu$. A {\bf $\lmu$-forest of rank $r$} 
is a~marked forest of rank $r$, $(T,\eta)$ such that $\mu=\lam_{r+1}(T,\eta)$ 
and $\lam\vdash\lam_0(T,\eta)$. 
\end{df}

  We call $((2,1^{n-2}),\mu)$-forests simply {\it $\mu$-forests} and 
$((2,1^{n-2}),(n))$-forests simply {\it $n$-trees}.
 
\vskip-10pt
\hskip-100pt $$\epsffile{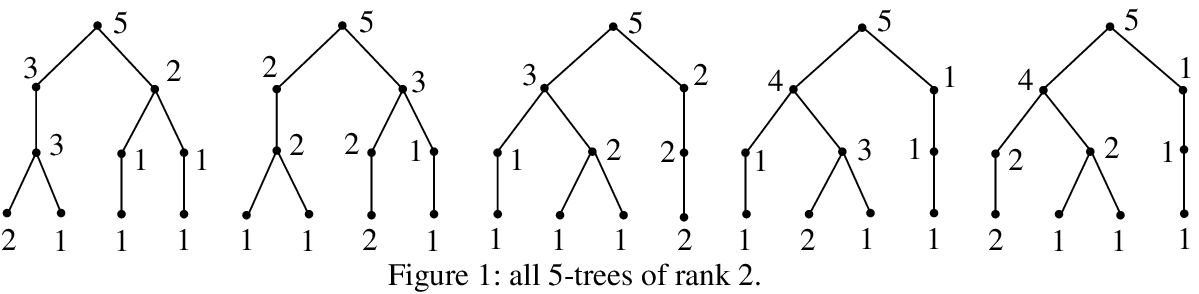}$$
\vskip-5pt
 
  Whenever $(T,\eta)$ is a~$\lmu$-forest of rank $r$ and $0\leq i\leq r$, we
can obtain a~$\lmu$-forest $(T^i,\eta^i)$ of rank $r-1$ by deleting from $T$
all the vertices of height~$i$ and connecting the vertices of height~$i+1$ to 
their grandchildren (unless $i=0$); $\eta^i$ is the restriction of 
$\eta$ to $V(T^i)$. In other words, $(T^i,\eta^i)$ is obtained from $(T,\eta)$ 
by removing the entire $i$th level, counting from the leaves, and filling in 
the gap in an~obvious way. This allows us to define a~boundary operator by
\begin{equation}\label{eqb}
  \partial(T,\eta)=\sum_{i=0}^r(-1)^i(T^i,\eta^i).
\end{equation}
For example:

\vskip-10pt
\hskip-100pt $$\epsffile{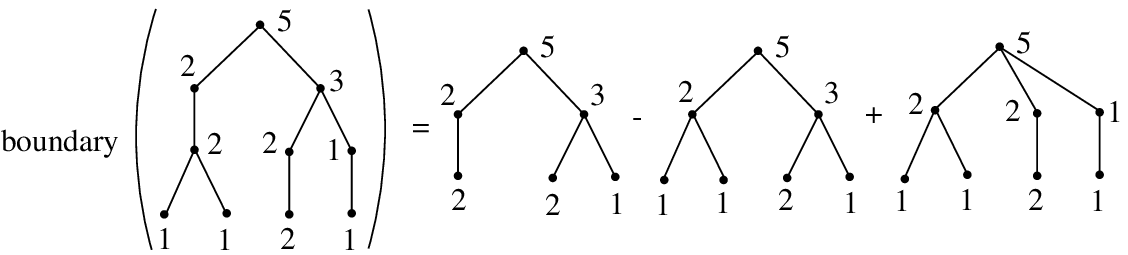}$$
\vskip-5pt

  For a~given set partition $\pi$ one can define the notion of 
a~{\bf $\pi$-forest $(T,\zeta)$ of rank $r$} almost identically to 
the case of number partitions described above. The difference is that 
$\zeta$ maps $V(T)$ to the set of finite sets, rather than the set of 
natural numbers. The condition~\eqref{mc1} is replaced~by 
\begin{equation}\label{mc2}
\zeta(v)=\bigcup_{w\in\text{\rm{children}}(v)}\zeta(w),
\end{equation}
and $\pi=\{\zeta(v)\,|\,v\in R(T)\}$. For $0\leq i\leq r+1$, analogously to 
$\lam_i(T,\eta)$, we define $\pi_i(T,\zeta)$ to be the set partition which is read
off from the vertices of $T$ having height~$i$. 

  Let $\mu$ be the type of $\pi$, then there exists a~canonical $\mu$-forest 
$(T,|\zeta|)$ associated to each $\pi$-forest $(T,\zeta)$, where $|\zeta|$ 
is obtained as the composition of $\zeta$ with the map which maps finite sets to 
their sizes.

\subsection{The main theorem.}$\,$

\noindent
  Let us describe how to associate a $\lmu$-forest, $\psi(\sigma)$, of rank $r$ 
to an $r$-simplex $\sigma$ of $\xlm$. The simplex $\sigma$ is an~$\st_\pi$-orbit 
of $r$-simplices of $\da(\Pi_\lam(\hat 0,\pi))$, where $\pi$ is a set partition 
of type~$\mu$. Take a representative of this orbit, a~chain $c=(x_r>\dots>x_0)$. 
Now we define $\psi(\sigma)=(T,\eta)$. Each element $x_i$ corresponds to the
$i$th level in $T$, counting from the leaves. Each block $b$ of $x_i$ 
corresponds to a node in the tree; on this node we define the value of~$\eta$ 
to be $|b|$. We define the edges of the tree $T$ by connecting each node 
corresponding to a block $b$ of $x_i$ to all nodes corresponding to the 
blocks of $x_{i-1}$ contained in $b$, we do that for all $b$ and $i$.
The top $(r+1)$th level is added artificially, its nodes correspond to the 
blocks of $\pi$, and the egdes from the top level to the $r$th level 
connect each block of $\pi$ to the blocks of $x_r$ contained in it. 
For example, the value of $\psi$ on the $\cs_5$-orbit of the chain 
$(123)(45)>(123)(4)(5)>(13)(2)(4)(5)$ is the first 5-tree on the Figure~1.

  We are now ready to state and prove the main result of this section.

\begin{thm}\label{main1}
  Assume $\lam\vdash\mu\vdash n$, $\lam\neq\mu$. The correspondence $\psi$ of 
the $r$-simplices of $\xlm$ and $\lmu$-forests $(T,\eta)$ of rank $r$ is 
a~bijection. Under this bijection, the boundary operator of the triangulated 
space $\xlm$ corresponds to the boundary ope\-ra\-tor described in~\eqref{eqb}.
\end{thm}
  
   In particular, the simplices of $\da(\Pi_n)/\cs_n$ along with the cell inclusion
structure are described by the $n$-trees. Indeed, $\da(\Pi_5)/\cs_5$ is shown in the
figure below.
  
\vskip-5pt
\hskip-100pt $$\epsffile{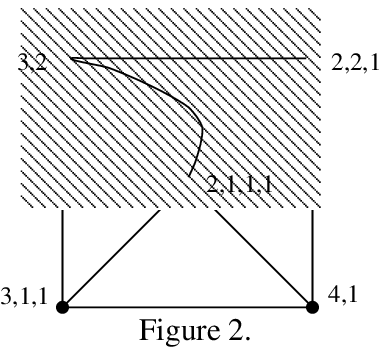}$$
\vskip-5pt

 The five triangles may be labeled by the five 5-trees of rank 2 in Figure 1.

\noindent
{\bf Proof of the Theorem \ref{main1}.} By the definitions of $\Pi_n$ and of $\da$, 
the $r$-simplices of $\da(\Pi_n)$ can be indexed by $([n])$-trees of rank~$r$
(we write $([n])$ to emphasize that the set $[n]$ is viewed here as a~set partition
consisting of only one set). Furthermore, the cell inclusions in $\da(\Pi_n)$
correspond to level deletion in $([n])$-trees as is described above for the case
of number partitions, because the levels in the $([n])$-trees correspond to the
elements of $\Pi_n$, and the edges in the $([n])$-trees correspond to block
inclusions of two consecutive elements in the chain. 

  More generally, the $r$-simplices $\sig$ of $\da(\pl(\hat 0,\pi))$ can be 
indexed by $\pi$-forests $(T(\sig),\zeta(\sig))$ such that $\lam$ refines 
the type of $\pi_{0}(T(\sig),\zeta(\sig))$. The definition of $\psi$ can 
now be rephrased as {\it associating to $\sigma$ the $\lmu$-forest 
$(T(\sig),|\zeta(\sig)|)$,} where $\mu=\,$type$\,\pi$.

  The group action of $\st_\pi$ on $\da(\pl(\hat 0,\pi))$ corresponds 
to relabeling elements within the sets of $\zeta(\sig)$. This shows that 
for $g\in\st_\pi$ we have $T(g\sig)=T(\sig)$ and $|\zeta(g\sig)|=|\zeta(\sig)|$.
Therefore $\psi(\sig)$ is well-defined, it does not depend on 
the choice of the representative of the, corresponding to $\sig$, 
$\st_\pi$-orbit of chains.

  $\psi$ is surjective, we shall now show that it is also injective. 
If $\sig_1,\sig_2$ are two different $r$-simplices of
$\da(\pl(\hat 0,\pi))$ such that $T(\sig_1)=T(\sig_2)$ and $|\zeta(\sig_1)|=
|\zeta(\sig_2)|$, then there exists $g\in\st_\pi$ such that $\zeta(g\sig_2)=
\zeta(\sig_1)$. Indeed, let $T=T(\sig_1)=T(\sig_2)$ and let $\alpha_1$, 
resp.~$\alpha_2$, be the string concatenated from the values of $\zeta(\sig_1)$,
resp.~$\zeta(\sig_2)$, on the leaves of $T$; the order of leaves of $T$ is arbitrary,
but the same for $T(\sig_1)$ and $T(\sig_2)$, the order of elements within each
$\zeta(\sig_1)(v)$, resp.~$\zeta(\sig_2)(v)$, for $v\in L(T)$ is also chosen
arbitrarily. Then $g\in\cs_n$ which maps $\alpha_2$ to $\alpha_1$ satisfies
the necessary conditions: $\zeta(g\sig_2)=\zeta(\sig_1)$ on the leaves
of~$T$, and hence by~\eqref{mc2} on all vertices of~$T$. Furthermore, since
$g\pi=g\pi_{r+1}(T,\zeta(\sig_2))=\pi_{r+1}(T,\zeta(g\sig_2))=
\pi_{r+1}(T,\zeta(\sig_1))=\pi$, we have $g\in\st_\pi$.

  This shows that $\psi$ is a~bijection. Since the levels of 
the $(\lam,\mu)$-forests correspond to the $\st_\pi$-orbits of 
the vertices of $\da(\Pi_\lam(\hat 0,\pi))$ (hence to the 
vertices of $\xlm$), the boundary operator of $\xlm$ corresponds
under $\psi$ to the level deletion in $\lmu$-forests, i.e.~the
boundary operator described in \eqref{eqb}.
\qed
 
\subsection{Remarks.} $\,$

1. While the presence of the root in an~$([n])$-tree is just a~formality
(two marked $([n])$-trees are equal iff the deletion of the root gives equal
marked forests), the presence of the roots in a~$\pi$-forest is vital. In fact,
if roots were not taken into account (as seems natural, since the partition read
off from the roots does not correspond to any vertex in $\da(\pl(\hat 0,\pi))$)
the argument above would be false already for vertices: if $\tau_1,\tau_2\in
\pl(\hat 0,\pi)$, such that type$\,(\tau_1)=\,$type$\,(\tau_2)$ (i.e.,~the 
corresponding marked forests of rank~0 are equal once the roots are removed),
there may not exist $g\in\st_\pi$, such that $g\tau_2=\tau_1$ (although such
$g\in\cs_n$ certainly exists).                
     
2. Marked forests equipped with an~order on the children of each vertex were used 
by Vassiliev, \cite{Vas}, to label cells in a~certain CW-complex structure on the 
space $\wti\dr^n(m)$, the one-point compactification of the configuration space 
of~$m$ unmarked distinct points in $\dr^n$. Vassiliev's cell decomposition of 
$\wti\dr^n(m)$ is a~generalization of the earlier Fuchs' cell decomposition of 
$\wti\dr^2(m)=\wti\dc(m)$, \cite{Fu}, which allowed Fuchs to compute the ring 
$H^*(Br(m),\dz_2)$, where $Br(m)$ is Artin's braid group on $m$ strings, see 
also~\cite{Vai}. Beyond a~certain similarity of the combinatorial objects used 
for labeling the cells, cf.~\cite[Lemma 3.3.1]{Vas} and Theorem~\ref{main1}, the 
connection between this paper and the results of Vassiliev and Fuchs seems unclear.

  As yet another instance of a~similar situation, we would like to mention the 
labeling of the components in the stratification of $\overline{M}_{0,n}$ (the 
Deligne-Knudsen-Mumford compactification of the moduli space of stable projective
complex curves of genus $0$ with $n$ punctures) with trees with $n$ labeled leaves,
see \cite{FM,Kn}.

           \section{A new proof of a~theorem of Arnold}   \label{s5}
           
\subsection{Formulation of the main theorem and corollaries.}$\,$

\noindent
   In this section we take a~look at a~rather general question of which
$\dq$-acyclicity of the spaces $\xlm$ is a~special case:
\begin{quote}
   Let $\pi\in\Pi_n$ and let $Q$ be an~$\st_\pi$-invariant subposet of
$\Pi_n(\hat 0,\pi)$. When is the multiplicity of the trivial representation
in the induced representation of $\st_\pi$ on $\wti H_i(\da(Q),\dq)$ equal 
to~0 for all $i$, in other words, when is $\da(Q)/\st_\pi$ $\dq$-acyclic?
\end{quote}
 
\begin{df}
   Let $\Lambda$ be a~subset of the set of all number partitions of $n$ such that
$(1^n),(n)\notin\Lambda$. Define $\Pi_\Lambda$ to be the subposet of $\Pi_n$
consisting of all set partitions $\pi$ such that $(\text{\rm type }\pi)\in\Lambda$.    
\end{df}     

  Clearly, $\Pi_\Lambda$ is $\cs_n$-invariant and, more generally, 
$\Pi_\Lambda(\hat 0,\pi)$ is $\st_\pi$-invariant. Vice versa, any $\cs_n$-invariant
subposet of $\Pi_n\sm\{(\{1\},\dots,\{n\}),([n])\}$ is of the form $\Pi_\Lambda$ for 
some~$\Lambda$.

  The following theorem is the main result of this section. It is a~generalization
of \cite[Theorem 4.1]{Ko2}  which now covers Theorem~\ref{arn}. The proof is
a~combination of the language of marked forests from Section~\ref{s3} and the ideas 
used in the proof of~\cite[Theorem~4.1]{Ko2}.

\begin{thm}\label{main2}
  Let $2\leq k<n$. Assume $\Lambda$ is a~subset of the set of all number 
par\-ti\-tions of~$n$ such that $(1^n),(n)\notin\Lambda$ and $\Lambda$ satisfies 
the following condition:
\begin{quote}
  {\rm Condition $C_k$.} If $\mu\in\Lambda$, such that $\mu=(\mu_1,\dots,\mu_t)$, 
where $\mu_i=kq_i+r_i$, $0\leq r_i<k$ for $i\in[t]$, then   
$\gamma_k(\mu)=(k^{q_1+\dots+q_t},1^{r_1+\dots+r_t})\in\Lambda$.  
\end{quote}
   Then for any $\mu\in\Lambda\cup\{(n)\}$ the triangulated space 
$\xLm=\da(\Pi_\Lambda(\hat 0,\pi))/\st_\pi$, where $\mu=\,${\rm type}$\,\pi$, is 
collapsible, in particular the multiplicity of the trivial representation in the 
induced $\st_\pi$-representation on $\wti H_i(\da(\Pi_\Lambda(\hat 0,\pi)),\dq)$ 
is equal to~0 for all~$i$.
\end{thm}

\begin{crl}
   Assume $\lam=(k^m,1^{n-km})$, $\lam\vdash\mu$, then $\xlm$ is collapsible.
In particular, $\xlm$ is $\dq$-acyclic, therefore Theorem~\ref{arn} follows.
\end{crl}
\pr Clearly $\xlm=\xLm$ for $\Lam=\{\tau\,|\,\lam\vdash\tau\vdash\mu,\tau\neq(n)\}$.
It is easy to check that Condition $C_k$ is satisfied for the case 
$\lam=(k^m,1^{n-km})$, therefore Theorem~\ref{arn} follows from Theorem~\ref{main2} 
via~\eqref{rswf}.
\qed

 Another consequence of Theorem~\ref{main2}, which was pointed out already 
in~\cite{Ko2}, has nothing to do with the spaces $\sla$. Namely, for $k=2$,
$\Lam=\{\tau\vdash n\,|\,r\leq l(\tau)\leq n-1\}$, for some $n-1\geq r\geq 2$, 
and $\mu=(n)$, we get results of Stanley ($r=2$),~\cite[page~151]{St}, and 
Hanlon ($r>2$),~\cite[Theorem~3.1]{Ha}.

  Theorem~\ref{main2} can be viewed as an~attempt to provide a~common framework for 
these results in the spirit of the question stated in the beginning of this section.
    
  \subsection{Auxiliary propositions.}$\,$

\noindent   
 First we need some terminology. For an arbitrary cell complex $\da$
we denote by $V(\da)$ the set of vertices of $\da$. Assume $\da$ is 
a regular CW complex and $\da'$ is its subcomplex. We denote the set 
of the simplices of $\da$ which are not simplices of $\da'$ by 
$\da\sm\da'$. We use the sign $\succ$ to denote the cover relation 
in the cell structure of~$\da$.

  Assume that, in addition, $\da$ is a triangulated space with some
linear order $\ll$ on the set of vertices. For $\sig\in\da\sm\da'$
we may write $\sig=(x_1,\dots,x_t)$, this notation is slightly
inaccurate since the set of vertices does not determine the simplex 
uniquely, all we mean is that $\sig$ has vertices $x_1\ll\dots\ll x_t$. 
In that case, we let $\xi(\sig)=i$ if $x_1,\dots,x_{i-1}\in V(\da')$ 
and $x_i\notin V(\da')$.

\begin{prop}\label{pr1}
   Let $\da$ be a~regular CW complex and $\da'$ a~subcomplex of $\da$, then the
   following are equivalent:
   
 a) there is a~sequence of collapses leading from $\da$ to $\da'$;
 
 b) there is a~matching of cells of $\da\sm\da'$: $\sig\leftrightarrow
    \phi(\sig)$, such that $\phi(\sig)\succ\sig$ and there is no sequence 
    $\sig_1,\dots\sig_t\in\da\sm\da'$ such that $\phi(\sig_1)\succ\sig_2,
    \phi(\sig_2)\succ\sig_3,\dots,\phi(\sig_t)\succ\sig_1$ (such matching is 
    called acyclic).
\end{prop}
\pr a) $\Rightarrow$ b). Let the elementary collapses define the matching $\phi$.
Assume there is a~sequence $\sig_1,\dots,\sig_t\in\da\sm\da'$ such that 
$\phi(\sig_1)\succ\sig_2,\phi(\sig_2)\succ\sig_3,\dots,\phi(\sig_t)\succ\sig_1$. 
Without loss of generality we can assume that the collapse $(\sig_1,\phi(\sig_1))$ 
precedes collapses $(\sig_i,\phi(\sig_i))$ for $2\leq i\leq t$.
Then $\phi(\sig_t)\succ\sig_1$ yields a~contradiction.

b) $\Rightarrow$ a). The proof is again very easy, various versions of it were
given in \cite[Corollary 3.5]{Fo}, \cite[Proposition 3.7]{Betc}, and 
\cite[Theorem 3.2]{Ko2}.
\qed  

\begin{prop}\label{pr2}
   Let $\da$ be a~triangulated space with some linear order $\ll$ on its set of 
vertices $V(\da)$. Let $V'\subseteq V(\da)$ and $\da'$ be the subcomplex of $\da$ 
induced by~$V'$. Assume we have a~partition $V(\da)=\cup_{z\in V'}V_z$ 
such that $z=\min_{\ll}V_z$. For $\sig\in\da\sm\da'$, let 
$\chi(\sig)\in V'$ be defined by $x_{\xi(\sig)}\in V_{\chi(\sig)}$. 
Finally assume that the following condition is satisfied:
\begin{quote}
  Condition $\aleph$. If $\sig\in\da\sm\da'$, $\sig=(x_1,\dots,x_t)$, is such that 
either $\xi(\sig)=1$ or $x_{\xi(\sig)-1}\neq\chi(\sig)$, then there exists a~unique 
simplex $\sig'=(x_1,\dots,x_{\xi(\sig)-1},\chi(\sig),x_{\xi(\sig)},\dots,x_t)$ such 
that $\sig'\sm\chi(\sig)=\sig$.
\end{quote}
  Then there is a~sequence of collapses leading from $\da$ to $\da'$.
\end{prop}
\pr Let $U$ denote the set of all $\sig\in\da\sm\da'$, $\sig=(x_1,\dots,x_t)$, 
such that $x_{\xi(\sig)-1}\neq\chi(\sig)$ or $\xi(\sig)=1$. The matching $\phi$ is 
defined by Condition $\aleph$: for $\sig\in U$ we set $\phi(\sig)=\sig'$. By 
Proposition~\ref{pr1} it is enough to check that this matching is acyclic.

  For $\sig\in U$ we have $\xi(\phi(\sig))=\xi(\sig)+1$. Moreover, if 
$\phi(\sig)\succ\sig'$ and $\sig'\in U$, then $\sig'=\phi(\sig)\sm x_{\xi(\sig)}$, 
hence $\xi(\sig')\geq\xi(\phi(\sig))$. Therefore, if there is a~sequence 
$\sig_1,\dots,\sig_t\in\da\sm\da'$ such that $\phi(\sig_1)\succ\sig_2,\phi(\sig_2)
\succ\sig_3,\dots,\phi(\sig_t)\succ\sig_1$, then we have $\xi(\sig_1)<
\xi(\phi(\sig_1))\leq\xi(\sig_2)<\xi(\phi(\sig_2))\leq\dots<\xi(\phi(\sig_t))\leq
\xi(\sig_1)$ which yields a~contradiction.
\qed 
   
\subsection{Proof of the Theorem~\ref{main2}.}$\,$

\noindent
  We define a~$\Lmu$-forest of rank $r$ to be a~marked forest $(T,\eta)$ of rank $r$
such that $\lam_{r+1}(T,\eta)=\mu$ and $\lam_i(T,\eta)\in\Lam$, for $0\leq i\leq r$. 
It follows from the discussion in Section~\ref{s3} and in particular from 
Theorem~\ref{main1} that the $r$-simplices of $\xLm$ can be indexed by $\Lmu$-forests
of rank $r$ so that the boundary relation of $\xLm$ corresponds to level
deletion in the marked forests. 

  We call number partitions of the form $(k^m,1^{n-km})$, for some~$m$, 
{\it special}. Let $K$ be the subcomplex of $\xLm$ induced by the set 
of all special partitions. We adopt the notations $\xi(\sig)$ and $\chi(\sig)$ 
used in Proposition~\ref{pr2} to the context of $\xLm$ and its subcomplex $K$. 
The linear order on $V(\xLm)$ can be taken to be any linear extension of the
partial order on $V(\xLm)$ given by the negative of the length function.
The partition of $V(\xLm)$ is given by: for $v\in V(\xLm)\sm V(K)$, 
$z\in V(K)$, we have $v\in V_z$ iff $z=\gamma_k(v)$.

  Let us show that the subcomplex $K$ satisfies Condition $\aleph$. Let
$\sig\in\xLm\sm K$, $\sig=(x_1,\dots,x_t)$, and assume $\xi(\sig)=1$ or $\chi(\sig)
\neq x_{\xi(\sig)-1}$. In the language of marked forests this can be reformulated 
as: $\sig$ is a~$\Lmu$-forest $(T,\eta)$ of rank $t$ such that $\lam_{\xi(\sig)-1}
(T,\eta)$ is not special and if $\xi(\sig)>1$ then $\lam_0(T,\eta),\dots,
\lam_{\xi(\sig)-2}(T,\eta)$ are special, and $\lam_{\xi(\sig)-2}(T,\eta)\neq
\gamma_k(\lam_{\xi(\sig)-1}(T,\eta))$.
In other words, on all vertices of height 0 to $\xi(\sig)-2$ the function~$\eta$
takes only values 1 or $k$ and for the vertices of height $\xi(\sig)-1$ it is no 
longer true. Moreover, there exists a~vertex of height $\xi(\sig)-1$ which has at 
least~$k$ children on which $\eta$ is equal to~1. It is now clear that there exists 
a~{\it unique} $\Lmu$-forest $(\wti T,\ti\eta)$ of rank $r+1$ such that 
\begin{itemize}
\item $(\wti T^{\xi(\sig)-1},\ti\eta^{\xi(\sig)-1})=(T,\eta)$;
\item $\lam_{\xi(\sig)-1}(\wti T,\ti\eta)=\gamma_k(\lam_{\xi(\sig)}(\wti T,\ti\eta))$,
i.e.,~$\ti\eta$ takes only values~1 or~$k$ on the vertices of height $\xi(\sig)-1$ 
and each vertex of height $\xi(\sig)$ in $(\wti T,\ti\eta)$ has no more than~$k-1$
children labeled~1.
\end{itemize}
   To construct $(\wti T,\ti\eta)$, extend $(T,\eta)$ by splitting each vertex
of height $\xi(\sig)-1$ into vertices marked~$k$ and~1 so that the number of $k$'s
is maximized. The uniqueness of $(\wti T,\ti\eta)$ follows from the definition
of the notion of isomorphism of marked forests.
 
   We have precisely checked Condition $\aleph$ and therefore by 
Proposition~\ref{pr1} we conclude that there is a~sequence of collapses leading
from $\xLm$ to $K$.

  It remains to see that $K$ is collapsible. If $\mu=(n)$, then $K$ is
a~simplex, so we can assume that $\mu\in\Lam$. If $\mu=\gamma_k(\mu)$,
then $K$ is again a simplex. Otherwise it is easy to see that there is 
a~unique vertex in $\xLm$ labeled $\gamma_k(\mu)$ and that $K$ is 
a~cone with an~apex in this vertex.
\qed

              \section{On the conjecture of Sundaram and Welker}
           \label{s4}    

\subsection{A counterexample to the general conjecture.} $\,$
 
\noindent
 The original formulation of Conjecture~\ref{c3} in \cite{SuW} was
  \begin{conj}\label{c1}\cite[Conjectures 4.12 and 4.13]{SuW}.
 Let $\lam$ and $\mu$ be different set par\-ti\-tions, such that $\lam\vdash\mu$. 
Let $\pi\in\pl$ be a~par\-tition of type $\mu$. Then the multiplicity of the trivial 
representation in the $\st_\pi$-module $\wti H_*(\da(\pl(\hat 0,\pi)),\dq)$~is~$0$.
\end{conj}

 In our terms Conjecture \ref{c1} is equivalent to
\begin{conj}\label{c2}
  For $\lam\vdash\mu$, $\lam\neq\mu$, the space $\xlm$ is $\dq$-acyclic.
\end{conj}
               
  We shall give an~example
when $\xlm$ is not even connected. It turns out that if one is only interested
in counting the number of connected components of $\xlm$, then there is a~simpler
poset model which we now proceed to describe. 

\begin{df}\label{pm}
  Assume $\lam\vdash\mu\vdash n$, $\lam\neq\mu$. The $\lmu$-forests of rank $0$ can 
be partially ordered as follows: $(T_1,\eta_1)\prec(T_2,\eta_2)$ if there exists
a~$\lmu$-forest $(T,\eta)$ of rank $1$ such that $(T_1,\eta_1)=(T^1,\eta^1)$ and 
$(T_2,\eta_2)=(T^0,\eta^0)$. We call the obtained poset $P_{\lam,\mu}$.
\end{df}      
   In other words, elements of $P_{\lam,\mu}$ are number partitions $\tau\neq\mu$ 
such that $\lam\vdash\tau\vdash\mu$, together with a~bracketing which shows how 
to form $\mu$ out of $\tau$, the order of the brackets and of the terms within 
the brackets is neglected. For example $(1,1,1)(3,1)(2,2)$ and $(3)(2,1,1)(2,1,1)$ 
are two different elements of $P_{(2,1^9),(4^2,3)}$, while $(1,1,1)(2,2)(3,1)$ is 
equal to the first mentioned element. These bracketed partitions are ordered by 
refinement, preserving the bracket structure.

\begin{prop}\label{pr}
  $X_{\lam,\mu}$ and $\da(P_{\lam,\mu})$ have the same number of connected 
components, i.e., $\beta_0(X_{\lam,\mu})=\beta_0(\da(P_{\lam,\mu}))$.
\end{prop}

\pr We know that $\da(P_{\lam,\mu})$ and $X_{\lam,\mu}$ have the same set of vertices 
and that there is an~edge between two vertices $a$ and $b$ of $\da(P_{\lam,\mu})$ iff
$a\prec b$ or $b\prec a$, which is, by the Definition \ref{pm}, the case iff there is 
an~edge between the corresponding vertices of $X_{\lam,\mu}$. This shows that 
$\da(P_{\lam,\mu})$ and $\xlm$ have the same number of connected components.
\qed

   Note that 1-skeleta of $\xlm$ and $\da(P_{\lam,\mu})$ need not be equal. 
$\da(P_{\lam,\mu})$ can intuitively be thought of as a~simplicial complex obtained
by forgetting the multiplicities of simplices in the triangulated space $\xlm$. 

\noindent
{\bf Counterexample.} For $n=23$, $\lam=(7,6,4,3,2,1)$, $\mu=(10,8,5)$, $\xlm$ is 
disconnected. $P_{\lam,\mu}$ is shown on the figure below. Clearly $\da(P_{\lam,\mu})$ 
is not connected, hence, by the~Proposition~\ref{pr}, neither is $\xlm$, which 
disproves Conjecture~\ref{c1}.
 
\vskip-5pt
\hskip-100pt $$\epsffile{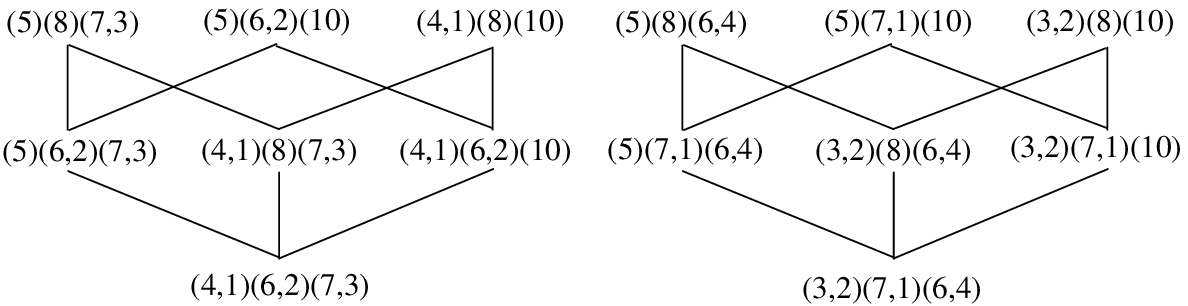}$$
\vskip-5pt

{\it Remark.} In the counterexample above, one can actually verify that 
$\xlm=\da(P_{\lam,\mu})$. However, we choose to use posets $P_{\lam,\mu}$ for
two reasons: 
\begin{enumerate}
\item it is easier to produce series of counterexamples to Sundaram-Welker
conjecture using $\da(P_{\lam,\mu})$ rather than $\xlm$; 
\item we feel that posets $P_{\lam,\mu}$ are of independent interest, since they 
are in a~certain sense the "naive" quotient of $\Pi_\lam(\hat 0,\pi)$ by $\st_\pi$.
\end{enumerate}

 We believe that, in general, connected components of $\xlm$ may be
 not acyclic.

\subsection{Proof of the conjecture in a special case.}
  
\begin{df}\label{dgen}
  We say that a~number partition $\lam=(\lam_1,\dots,\lam_t)$ is generic
(also called free of resonances in~\cite{ShW}, having no equal subsums in~\cite{Ko1}) 
if whenever $\sum_{i\in I}\lam_i=\sum_{j\in J}\lam_j$, for some $I,J\subseteq[t]$,
we have $\{\lam_i\}_{i\in I}=\{\lam_j\}_{j\in J}$ as multisets.
\end{df}

  For example $\lam=(k^m)$ is generic.
  
\begin{thm}\label{tgen}  
  If $\lam$ is generic, we have 
$$\tb_i(\sla,\dq)=\begin{cases} 1,& \text{for } i=2l(\lam);\\
                                0,& \text{otherwise.}
  \end{cases}$$
\end{thm}
\pr As we pointed out before, it is enough to show that $\xlm$ is $\dq$-acyclic
for $\lam\vdash\mu$, $\lam\neq\mu$. 

  If $\lam$ is generic, any number $1\leq m\leq n$ can be partitioned into numbers 
from $\lam$ in at most one way. First, it implies that there exists a~unique 
$\lmu$-forest $(T,\eta)$ of rank $0$ such that $\lam_0(T,\eta)=\lam$, denote it by~$x$. 
Second, for any $\lmu$-forest $(T,\eta)$ of rank~$r$ such that $\lam_0(T,\eta)\neq\lam$ 
there exists a~unique $\lmu$-forest $(\wti T,\ti\eta)$ of rank~$r+1$ such that 
$(\wti T^0,\ti\eta^0)=(T,\eta)$ and $\lam_0(\wti T,\ti\eta)=\lam$. 

  In terms of the triangulated space $\xlm$ this means that there exists
a~vertex~$x$ such that each $r$-simplex of $\xlm$ not containing $x$ is 
contained in a~unique $(r+1)$-simplex containing~$x$. This means that, 
whenever $\lam$ is generic and $\lam\vdash\mu$, $\lam\neq\mu$, $\xlm$ is 
a~cone, in particular it is $\dq$-acyclic.
\qed

\vskip2pt
\noindent
{\bf Acknowledgments.} We would like to thank Eva-Maria Feichtner and
the anonymous referee for the careful reading of this paper.

\end{document}